\documentclass{amsart}

\usepackage{amsmath,amsthm,amsfonts,amscd,amssymb}

\newtheorem{thm}{Theorem}[section]
\newtheorem{cor}[thm]{Corollary}
\newtheorem{lem}[thm]{Lemma}
\newtheorem{prop}[thm]{Proposition}

\newtheorem{prob}{Problem}

\theoremstyle{definition}

\theoremstyle{remark}

\begin{document}

\title{Integral Points of Small Height Outside of a Hypersurface}
\author{Lenny Fukshansky}

\address{Department of Mathematics, Texas A\&M University, College Station, Texas 77843-3368}
\email{lenny@math.tamu.edu}
\subjclass{Primary 11C08, 11H06; Secondary 11D04, 11H46}
\keywords{polynomials, lattices, linear forms, height}

\begin{abstract}
Let $F$ be a non-zero polynomial with integer coefficients in $N$ variables of degree $M$. We prove the existence of an integral point of small height at which $F$ does not vanish. Our basic bound depends on $N$ and $M$ only. We separately investigate the case when $F$ is decomposable into a product of linear forms, and provide a more sophisticated bound. We also relate this problem to a certain extension of Siegel's Lemma as well as to Faltings' version of it. Finally we exhibit an application of our results to a discrete version of the Tarski plank problem.
\end{abstract}

\maketitle

\def\A{{\mathcal A}}
\def\B{{\mathcal B}}
\def\C{{\mathcal C}}
\def\D{{\mathcal D}}
\def\F{{\mathcal F}}
\def\x{{\mathcal H}}
\def\I{{\mathcal I}}
\def\J{{\mathcal J}}
\def\K{{\mathcal K}}
\def\L{{\mathcal L}}
\def\M{{\mathcal M}}
\def\R{{\mathcal R}}
\def\s{{\mathcal S}}
\def\V{{\mathcal V}}
\def\X{{\mathcal X}}
\def\Y{{\mathcal Y}}
\def\H{{\mathcal H}}
\def\cee{{\mathbb C}}
\def\pee{{\mathbb P}}
\def\que{{\mathbb Q}}
\def\real{{\mathbb R}}
\def\zed{{\mathbb Z}}
\def\qbar{{\overline{\mathbb Q}}}
\def\eps{{\varepsilon}}
\def\ahat{{\hat \alpha}}
\def\bhat{{\hat \beta}}
\def\gt{{\tilde \gamma}}
\def\h{{\tfrac12}}
\def\be{{\boldsymbol e}}
\def\bei{{\boldsymbol e_i}}
\def\bc{{\boldsymbol c}}
\def\bm{{\boldsymbol m}}
\def\bk{{\boldsymbol k}}
\def\bi{{\boldsymbol i}}
\def\bl{{\boldsymbol l}}
\def\bq{{\boldsymbol q}}
\def\bu{{\boldsymbol u}}
\def\bt{{\boldsymbol t}}
\def\bs{{\boldsymbol s}}
\def\bv{{\boldsymbol v}}
\def\bw{{\boldsymbol w}}
\def\bx{{\boldsymbol x}}
\def\bX{{\boldsymbol X}}
\def\bz{{\boldsymbol z}}
\def\bwy{{\boldsymbol y}}
\def\bY{{\boldsymbol Y}}
\def\bL{{\boldsymbol L}}
\def\ba{{\boldsymbol\alpha}}
\def\bb{{\boldsymbol\beta}}
\def\bet{{\boldsymbol\eta}}
\def\bxi{{\boldsymbol\xi}}
\def\bo{{\boldkey 0}}
\def\bol{{\boldkey 1}_L}
\def\ep{\varepsilon}
\def\p{\boldsymbol\varphi}
\def\q{\boldsymbol\psi}
\def\rank{\operatorname{rank}}
\def\aut{\operatorname{Aut}}
\def\lcm{\operatorname{lcm}}
\def\sgn{\operatorname{sgn}}
\def\spn{\operatorname{span}}
\def\md{\operatorname{mod}}
\def\Norm{\operatorname{Norm}}
\def\dim{\operatorname{dim}}
\def\det{\operatorname{det}}
\def\Vol{\operatorname{Vol}}
\def\rk{\operatorname{rk}}

\section{Introduction and notation}

Let
$$F(\bX) = F(X_1,...,X_N) \in \zed[X_1,...,X_N]$$
be a polynomial in $N \geq 1$ variables of degree $M \geq 1$ with integer coefficients. If $F$ is not identically zero, there must exist a point with integer coordinates at which $F$ does not vanish, in other words an integral point that lies outside of the hypersurface defined by $F$ over $\que$. How does one find such a point?
\smallskip

For a point $\bx = (x_1,...,x_N) \in \zed^N$, define its {\it height} and {\it length} respectively by
$$H(\bx) = \max_{1 \leq i \leq N} |x_i|,\ \ \L(\bx) = \sum_{i=1}^N |x_i|.$$
It is easy to see that a set of points with height or length bounded by some fixed constant is finite. In fact, for a positive real number $R$,
\begin{equation}
\label{1.1}
\left| \{ \bx \in \zed^N : H(\bx) \leq R \} \right| = (2[R]+1)^N,
\end{equation}
and
\begin{equation}
\label{1.2}
\left| \{ \bx \in \zed^N : \L(\bx) \leq R \} \right| = \sum_{k=0}^{\min([R],N)} 2^k \binom{N}{k} \binom{[R]}{k},
\end{equation}
where (\ref{1.1}) is obvious and (\ref{1.2}) follows from Theorem 6 of \cite{local:riemann}; we write $[R]$ for the integer part of $R$. Therefore if we were able to prove the existence of a point $\bx \in \zed^N$ with $H(\bx) \leq R$ or $\L(\bx) \leq R$ for some explicitly determined value of $R$, then the problem of finding this point would reduce to a finite search. Thus we will consider the following problem.

\begin{prob} \label{problem1} Given a polynomial $F(X_1,...,X_N) \in \zed[X_1,...,X_N]$ in $N \geq 1$ variables of degree $M \geq 1$, prove the existence of a point $\bx \in \zed^N$ of small height (length) such that $F(\bx) \neq 0$ with an explicit upper bound on the height (length).
\end{prob}

This problem, although of independent interest, is also closely related to an important trend in Diophantine approximations. Consider the special case when $F(\bX)$ is homogeneous in $N \geq 2$ variables and decomposable into a product of integral linear forms, i.e. suppose that
$$F(X_1,...,X_N) = \prod_{i=1}^M L_i(X_1,...,X_M),$$
where $L_i(X_1,...,X_N) \in \zed[X_1,...,X_N]$ is a linear form for each $1 \leq i \leq M$. Write
$$\Lambda_i = \{ \bwy \in \zed^N : L_i(\bwy) = 0\},$$
for each $1 \leq i \leq M$. Then each such $\Lambda_i$ is a lattice of rank $N-1$ in $\real^N$, and the statement that $F(\bx) \neq 0$ for some $\bx \in \zed^N$ is equivalent to the statement that $\bx \in \zed^N \setminus \bigcup_{i=1}^M \Lambda_i$. In this case Problem \ref{problem1} can be restated as follows.

\begin{prob} \label{problem2} Let $\Lambda_1,...,\Lambda_M$ be sublattices of $\zed^N$ of rank $N-1$. Prove the existence of a point $\bx \in \zed^N \setminus \bigcup_{i=1}^M \Lambda_i$ of small height with an explicit bound on height.
\end{prob}

More generally, suppose that $\Omega \subseteq \zed^N$ is a lattice of rank $J$, $2 \leq J \leq N$, and let $\Lambda_1,...,\Lambda_M$ be proper sublattices of $\Omega$ of respective ranks $l_1,...,l_M$, $1 \leq l_i \leq J-1$ for each $1 \leq i \leq M$. We can state the following generalization of Problem \ref{problem2}.

\begin{prob} \label{problem3} Prove the existence of a point $\bx \in \Omega \setminus \bigcup_{i=1}^M \Lambda_i$ of small height with an explicit bound on the height.
\end{prob}

Notice that this can be viewed as a certain extension of a classical Siegel's Lemma (see for instance \cite{vaaler:siegel}, \cite{siegel}, and \cite{thue}), which guarantees the existence of a non-zero point of small height in $\Omega$. Then Problem \ref{problem3} is a version of Siegel's Lemma with additional linear conditions, and Problem \ref{problem2} can be thought of as a problem inverse to Siegel's Lemma. Here is the main result of this paper.

\begin{thm} \label{classical:subspace} Let $\Omega \subseteq \zed^N$ be a sublattice of rank $J$, $1 \leq J \leq N$. Let $\Lambda_1,...,\Lambda_M$ be proper non-zero sublattices of $\Omega$. Then there exists $\bx \in \Omega \setminus \bigcup_{i=1}^M \Lambda_i$ such that
\begin{equation}
\label{class_subspace:bound}
H(\bx) \leq \left( \left( \frac{3}{2} \right)^{J-1} J^J \left\{ \sum_{i=1}^M \frac{1}{H(\Lambda_i)} + \sqrt{M} \right\} + J \right) H(\Omega).
\end{equation}
\end{thm}

A generalization of this theorem to the number field case is presented in the companion paper \cite{me:number}. The paper is organized as follows. In section 2 we treat Problem \ref{problem1} for a general polynomial in $N$ variables of degree $M$ with integer coefficients and produce some basic results. The bounds on height and length depend on $N$ and $M$ only. In fact, we produce simple bounds on both, height and length, that depend only on $M$. In sections 3 and 4 we prove some technical lemmas on counting lattice points in simple convex bodies in a Euclidean space. These lemmas are then used in section 5 to treat Problem \ref{problem2} and Problem \ref{problem3}, proving Theorem \ref{classical:subspace}. The bounds on height of the points in question in section 5 depend on $M$, $N$, and heights of corresponding lattices. We also derive a sharper bound for Faltings' version of Siegel's Lemma (see Proposition 2.18, \cite{faltings}) in our context as a corollary of the main result. Finally in section 6 we exhibit a simple application of our results to a discrete version of the Tarski plank problem.
\bigskip

Before we can proceed we need to introduce some additional notation. We can extend height to polynomials by evaluating height of their coefficient vectors. We also define a height function on lattices. Let $\Lambda \subseteq \zed^N$ be a lattice of rank $J$, $1 \leq J \leq N$. Choose a basis $\bx_1,...,\bx_J$ for $\Lambda$, and write $X = (\bx_1\ ...\ \bx_J)$ for the corresponding $N \times J$ basis matrix. Then 
$$\Lambda = \{ X \bt : \bt \in \zed^J \}.$$
On the other hand, there exists an $(N-J) \times N$ matrix $A$ with integer relatively prime entries so that 
$$\Lambda \subseteq \{ \bx \in \zed^N : A \bx = 0 \}.$$
Let $\I$ be the collection of all subsets $I$ of $\{1,...,N\}$ of cardinality $J$. For each $I \in \I$ let $I'$ be its complement, i.e. $I' = \{1,...,N\} \setminus I$, and let $\I' = \{ I' : I \in \I\}$. Then 
$$|\I| = \binom{N}{J} = \binom{N}{N-J} = |\I'|.$$
For each $I \in \I$, write $X_I$ for the $J \times J$ submatrix of $X$ consisting of all those rows of $X$ which are indexed by $I$, and $_{I'} A$ for the $(N-J) \times (N-J)$ submatrix of $A$ consisting of all those columns of $A$ which are indexed by $I'$. By the duality principle of Brill-Gordan \cite{gordan:1} (also see Theorem 1 on p. 294 of \cite{hodge:pedoe})
\begin{equation}
\label{duality}
|\det (X_I)| = \gamma |\det (_{I'} A)|,
\end{equation}
for an appropriate $\gamma \in \que$. Then define the vectors of {\it Grassmann coordinates} of $X$ and $A$ respectively to be 
$$Gr(X) = (\det (X_I))_{I \in \I} \in \zed^{|I|},\ \ Gr(A) = (\det (_{I'} A))_{I' \in \I'} \in \zed^{|I'|},$$
and so by (\ref{duality})
$$H(Gr(X)) = \gamma H(Gr(A)).$$
Define height of $\Lambda$ denoted by $H(\Lambda)$ to be $H(Gr(X))$. This definition is legitimate, since it does not depend on the choice of the basis for $\Lambda$. In particular, notice that if $N \geq 2$ and
$$L(X_1,...,X_N) = \sum_{i=1}^N q_i X_i \in \zed[X_1,...,X_N]$$
is a linear form with a non-zero coefficient vector $\bq$ having relatively prime coordinates, and $\Lambda = \{ \bx \in \zed^N : L(\bx) = 0 \}$ is a sublattice of $\zed^N$ of rank $N-1$, then
\begin{equation}
\label{1.4}
H(\Lambda) = H(L) = H(\bq).
\end{equation}
We are now ready to proceed. Results of this paper also appear as a part of \cite{me:diss}. 
\bigskip

\section{Basic bounds}

In this section we prove the existence of integral points of small height and length at which a given polynomial in $N$ variables of degree $M$ does not vanish. Our bounds on height depend on $M$ and $N$ only. An application of results of this section is exhibited in \cite{me:smallzeros}.
\smallskip

First let
$$\M' = \M'(N,M) = \{ \bm \in \zed_{\geq 0}^N : m_1+...+m_N \leq M \}.$$
Then let
$$F(X_1,...,X_N) = \sum_{\bm \in \M'} f_{\bm} X_1^{m_1}...X_N^{m_N} \in \zed[X_1,...,X_N],$$
be a polynomial (not necessarily homogeneous) in $N \geq 1$ variables of degree $M \geq 1$ with coefficients in $\zed$. We write $\deg_{X_i} (F)$ for degree of $F$ in the variable $X_i$ for each $1 \leq i \leq N$, and $\deg(F)$ for the total degree of $F$. Let 
$$m(F) = \max_{1 \leq i \leq N} \deg_{X_i} (F),$$
then $m(F) \leq \deg(F) = M$.

\begin{lem} \label{2.2.1} Suppose $F(\bX)$ is not identically $0$. Let $S \subset \cee$ be a finite set of cardinality at least $m(F)+1$. Then there exists $\bx \in S^N$ such that $F(\bx) \neq 0$.
\end{lem}

This is an immediate corollary of Lemma 1 on p. 261 of \cite{cass:geom}. Notice that the assertion of Lemma \ref{2.2.1} is best possible (i.e. $|S|$ {\it must} be at least $m(F)+1$) as seen on the following example. Let $S = \{\alpha_1,...,\alpha_M\} \subset \zed$, and let
$$F(X_1,...,X_N) = \sum_{i=1}^N \prod_{j=1}^M (X_i - \alpha_j).$$
Then for each $\bx \in S^N$ we have $F(\bx) = 0$, where $m(F) = M = |S|$.

\begin{lem} \label{2.2.2} Let $F$ be as in Lemma \ref{2.2.1}. There exists $\bx \in \zed^N$ with $\bx_i \neq 0$ for all $1 \leq i \leq N$, $F(\bx) \neq 0$, and
\begin{equation}
\label{bound:M}
H(\bx) \leq \frac{M+2}{2}.
\end{equation}
Moreover, if $F$ is homogeneous, then the upper bound of (\ref{bound:M}) can be replaced with $\frac{M+1}{2}$.
\end{lem}

\proof
Recall that $M \geq m(F)$, and let
$$S = \left\{ -\left[ \frac{M}{2} \right]-1,...,-1,1,...,\left[ \frac{M}{2} \right]+1 \right\},$$
then $|S| = 2 \left( \left[ \frac{M}{2} \right]+1 \right) \geq M+1$. Hence, by Lemma \ref{2.2.1}, there must exist $\bx \in S^N$ such that $F(\bx) \neq 0$.
\smallskip

Now assume that $F$ is homogeneous and $M \geq 1$. Notice that if for any $1 \leq i \leq N$ the ``diagonal'' coefficient $f_{M \be_i} \neq 0$, then $F(\be_i) = f_{M \be_i} \neq 0$, and we are done. Hence assume $f_{M \be_i} = 0$ for all $1 \leq i \leq N$. Then each monomial of $F$ has degree $M$ and is a product of powers of at least two variables. Therefore $m(F) \leq M-1$, and so we can take
$$S = \left\{ -\left[ \frac{M-1}{2} \right]-1,...,-1,1,...,\left[ \frac{M-1}{2} \right]+1 \right\},$$
then $|S| = 2 \left( \left[ \frac{M-1}{2} \right]+1 \right) \geq M \geq m(F)+1$. Hence, by Lemma \ref{2.2.1}, there must exist $\bx \in S^N$ such that $F(\bx) \neq 0$. This completes the proof.
\endproof
\smallskip

Next we want to produce a bound on length of an integral point at which $F$ does not vanish. Consider $\bx$ of Lemma \ref{2.2.2}. Notice that 
\begin{equation}
\label{length:1}
\L(\bx) \leq N H(\bx) \leq \frac{N(M+2)}{2}.
\end{equation}
This is a trivial bound. We will produce a non-trivial bound on $\L(\bx)$ in a slightly more restrictive situation. Let $N \geq 2$, $M \geq 1$ be integers, and write
$$\M = \M(N,M) = \left\{ \bm \in \zed_{+}^N : m_1 + ... + m_N = M \right\}.$$
For the rest of this section, let
$$F(X_1,...,X_N) = \sum_{\bm \in \M} f_{\bm} X_1^{m_1}...X^{m_N} \in \zed[X_1,...,X_N],$$
be a non-zero homogeneous polynomial in $N$ variables of degree $M$ with coefficients in $\zed$.

\begin{lem} \label{2.2.4} Let $F$ be as above. There exists a point $\bx \in \zed^N$ such that $F(\bx) \neq 0$, and $\L(\bx) \leq \frac{(M+2)^2}{8}$.
\end{lem}

\proof
If $M=1$, then $F$ is just a linear form in $N \geq 2$ variables. Its nullspace has dimension $N-1$, and so cannot contain all the standard basis vectors. Therefore there exists $\bx \in \zed^N$ with $\L(\bx) = 1$ and $F(\bx) \neq 0$. From now on assume that $M \geq 2$. We can also assume that for each $1 \leq i \leq N$, the coefficient $f_{M\be_i} = 0$, where $\be_1,...,\be_N$ are the standard basis vectors, since if for some $1 \leq i \leq N$, $f_{M\be_i} \neq 0$, then $F(\be_i) = f_{M\be_i} \neq 0$.
\smallskip

We argue by induction on $N$. First suppose that $N=2$, then we can write
$$F(X_1,X_2) = \sum_{i=1}^{M-1} f_i X_1^i X_2^{M-i} = X_1 X_2 \sum_{i=1}^{M-1} f_i X_1^{i-1} X_2^{M-i-1},$$
and so $\frac{1}{X_1} F(X_1,1) = \sum_{i=1}^{M-1} f_i X_1^{i-1}$ is a polynomial in one variable of degree at most $M-2$, therefore it can have at most $M-2$ nonzero roots, and so there must exist an integer $\beta$ with $|\beta| \leq \frac{M-2}{2} + 1$ such that $F(\beta,1) \neq 0$. Then $\bx = (\beta,1)$ is the required point with 
$$\L(\bx) \leq \frac{M+2}{2} \leq \frac{(M+2)^2}{8},$$
since $M \geq 2$.
\smallskip

Next assume $N>2$. For each $1 \leq i \leq N$, define $F_i$, {\it $i$-th section of $F$}, to be the homogeneous polynomial in $N-1$ variables of degree $M$ obtained from $F$ by setting $i$-th variable equal to $0$. First suppose that all sections of $F$ are identically zero, then
$$F(X_1,...,X_N) = X_1...X_N G(X_1,...,X_N),$$
where $G$ is a homogeneous polynomial of degree $M-N$ (this is only possible if $N \leq M$). By Lemma \ref{2.2.2} and (\ref{length:1}), there exists $\bx \in \zed^N$ such that $x_j \neq 0$ for all $1 \leq j \leq N$, $G(\bx) \neq 0$, and
$$\L(\bx) \leq \frac{N}{2} (M-N+2).$$
Then $F(\bx) \neq 0$. Let
$$f(z) = \frac{z}{2} (M-z+2) = -\frac{1}{2} z^2 + \frac{(M+2)}{2} z,$$
then $f$ achieves its maximum when $z = \frac{M+2}{2}$, and this maximum value is $\frac{(M+2)^2}{8}$. Hence
\begin{equation}
\L(\bx) \leq \frac{(M+2)^2}{8}.
\end{equation}
Next assume that for some $1 \leq i \leq N$, $F_i$ is not identically zero. Then we are done by the inductive hypothesis. This completes the proof. 
\endproof
\bigskip

Notice that observations of (\ref{length:1}) and Lemma \ref{2.2.4} can be summarized as follows.

\begin{prop} \label{2.2.5} Let $F$ be as above. There exists a point $\bx \in \zed^N$ such that $F(\bx) \neq 0$, and 
\begin{equation}
\label{length:bound1}
\L(\bx) \leq \left[ \frac{(M+2)}{2}\ \min \left\{ N, \frac{M+2}{4} \right\} \right],
\end{equation}
where $[\ ]$, as above, stands for the integer part function.
\end{prop}
\bigskip

Notice that if $F$ is irreducible, then the bound of Proposition \ref{2.2.5} can be trivially improved by replacing $N$ with $N-1$ in the upper bound, since we can set any one variable equal to zero without making an irreducible polynomial identically zero. 
\bigskip

\section{Lattice points in an aligned box}

In the next two sections we produce estimates for the number of points of a sublattice of the integer lattice in a closed cube in $\real^N$. These estimates are later used to prove our main theorem.
\smallskip

Let $A = (a_{mn})$ be an $N \times N$, uppertriangular, nonsingular matrix with real entries. Let $u_m < v_m$ for $m=1, 2, \dots , N$ and write
$$B(\bu,\bv) = \{\bx\in\real^N: u_m < x_m \le v_m\}.$$
We will be interested in estimating the number of points $\bxi$ in $\zed^N$ such that $A\bxi$ belongs to the aligned box $B(\bu,\bv)$.  To
begin with we have the following special result.

\begin{lem} \label{2.3.1}  Assume that $a_{11} = a_{22} = \cdots = a_{NN} = 1$ and $v_m - u_m$ is a positive integer for each $m=1, 2, \dots , N$.
Then
\begin{equation}
\label{3.1}
|\{\bxi\in\zed^N: A\bxi\in B(\bu,\bv)\}| = \prod_{m=1}^N (v_m - u_m).
\end{equation}
\end{lem}

\proof
We argue by induction on $N$.  If $N = 1$ the result is trivial because one easily checks that the number of integer 
points $\xi_1$ such that $u_1 < \xi_1 \le v_1$ is equal to $[v_1] - [u_1]$.  As $v_1 - u_1$ is an integer we find that $[v_1] - [u_1] = v_1 - u_1$.

Now assume that $N\ge 2$.  Let $\bet$ be a point in $\zed^{N-1}$ with coordinates indexed by $n=2, 3, \dots , N$.  Then define
$$I_{N-1} = \{\bet\in\zed^{N-1}: u_m < \sum_{n=m}^N a_{mn}\eta_n \le v_m\ \text{for}\ m=2, 3, \dots , N\}.$$
By the inductive hypothesis we have
\begin{equation}
\label{3.2}
|I_{N-1}| = \prod_{m=2}^N (v_m - u_m).
\end{equation}
If $\bet$ is a point in $I_{N-1}$ then the number of integer points $\xi_1$ such that
\begin{equation}
\label{3.3}
u_1 < \xi_1 + a_{12}\eta_2 + a_{13}\eta_3 + \cdots + a_{1N}\eta_N \le v_1,
\end{equation}
is $v_1 - u_1$.  Clearly a point
$$\bxi = \begin{pmatrix}  \xi_1\\
                   \eta_2\\
                  \vdots\\
                   \eta_N\end{pmatrix}$$
satisfies the condition $A\bxi\in B(\bu,\bv)$ if and only if $\bet\in I_{N-1}$ and $\xi_1$ satisfies (\ref{3.3}).  We have shown that the number of 
such points is 
$$(v_1 - u_1)|I_{N-1}|,$$
and this proves the lemma.
\endproof

If we drop the condition that each edge length $v_m - u_m$ is an integer then we get the following estimates.

\begin{lem} \label{2.3.2}  Assume that $a_{11} = a_{22} = \cdots = a_{NN} = 1$, then
\begin{equation}
\label{3.4}
\prod_{m=1}^N [v_m - u_m] \le |\{\bxi\in\zed^N: A\bxi\in B(\bu,\bv)\}| \le \prod_{m=1}^N ([v_m - u_m] + 1).
\end{equation}
\end{lem}

\proof
When proving the lower bound on the left of (\ref{3.4}) we can assume that $1 \le v_m - u_m$ for each $m=1, 2, \dots , N$.  
Now select real numbers $u_m'$ and $v_m'$ so that
$$u_m \le u_m' < v_m' \le v_m\quad\text{and}\quad v_m' - u_m' = [v_m - u_m],\quad\text{for}\quad m=1, 2, \dots , N.$$
As $B(\bu',\bv') \subseteq B(\bu,\bv)$ the inequality follows from Lemma \ref{2.3.1}.  To obtain the upper bound on the right of (\ref{3.4}) we argue in
essentially the same way.  Select real numbers $u_m''$ and $v_m''$ so that 
$$u_m'' \le u_m < v_m \le v_m''\quad\text{and}\quad v_m'' - u_m'' = [v_m - u_m] + 1,\quad\text{for}\quad m=1, 2, \dots , N.$$
Then $B(\bu,\bv)\subseteq B(\bu'',\bv'')$, and again the inequality follows from Lemma \ref{2.3.1}.
\endproof

Next we drop the condition that the diagonal entries of the matirx $A$ are all equal to $1$.

\begin{cor} \label{2.3.3}  Assume that the diagonal entries $a_{11}, a_{22}, \dots , a_{NN}$ are all positive.  Then we have
\begin{equation}
\label{box:lattice}
\prod_{m=1}^N \Bigl[\frac{v_m - u_m}{a_{mm}}\Bigr] 
        \le |\{\bxi\in\zed^N: A\bxi\in B(\bu,\bv)\}| \le \prod_{m=1}^N \Bigl(\Bigl[\frac{v_m - u_m}{a_{mm}}\Bigr] + 1\Bigr).
\end{equation}
\end{cor}
\bigskip

\section{Lattice points in cubes}

In this section we focus on the case when the box of section 3 is actually a cube, and in this case extend the estimate of section 3 to lattices of not full rank.
\smallskip

\noindent
For the rest of this section, let $R \geq 1$, and define
$$C^N_R = \{\bx \in \real^N : \max_{1 \leq i \leq N} |x_i| \leq R \},$$
to be a cube in $\real^N$ centered at the origin with sidelength $2R$. Given a lattice $\Lambda \subseteq \zed^N$ of rank $N-l$ for $0 \leq l \leq N-1$ and determinant $\Delta$, we want to estimate the number of points of $\Lambda$ in $C_R^N$. First notice that if $\Lambda = A \zed^N$ has full rank, and matrix $A$ as in section 3 has fixed determinant $\Delta$, then the right hand side of (\ref{box:lattice}) takes its maximum value when $a_{mm} = 1$ for $N-1$ distinct values of $m$. This leads to the following corollary.

\begin{cor} \label{2.4.1} Let $\Lambda \subseteq \zed^N$ be a lattice of full rank in $\real^N$ of determinant $\Delta$. Then for each point $\bz$ in $\real^N$ we have
\begin{equation}
\label{cube:lattice}
|\Lambda \cap (C_R^N+\bz)| \leq \left(\frac{2R}{\Delta} + 1 \right) (2R + 1)^{N-1}.
\end{equation}
Moreover, if $R$ is a positive integer multiple of $\Delta$, we have
\begin{equation}
\label{cube:lattice1}
\frac{(2R)^N}{\Delta} \leq |\Lambda \cap (C_R^N+\bz)|.
\end{equation}
\end{cor}

Now suppose that $\Lambda$ of Corollary \ref{2.4.1} is not of full rank. The bounds we produce next, as well as the bound of Corollary \ref{2.4.1}, depend only on the lattice, not on the choice of the basis as in (\ref{box:lattice}).

\begin{thm} \label{2.4.2} Suppose that $\Lambda \subseteq \zed^N$ is a lattice of rank $N-l$, where $1 \leq l \leq N-1$. Let $\Delta$ be the maximum of absolute values of Grassmann coordinates of $\Lambda$, that is $\Delta=H(\Lambda)$. Then for each point $\bz$ in $\real^N$ we have
\begin{equation}
\label{notfull:rank}
|\Lambda \cap (C_R^N+\bz)| \leq \left( \frac{2R}{\Delta}+1 \right)(2R+1)^{N-l-1}.
\end{equation}
\end{thm}

\proof
Pick $\be_{k_1},...,\be_{k_l}$ such that the lattice 
$$\spn_{\zed} \{\Lambda,\be_{k_1},...,\be_{k_l} \} \subseteq \zed^N$$
has rank $N$. Write $X=(\bx_1 \hdots \bx_{N-l})$ for the $N \times (N-l)$ basis matrix of $\Lambda$. Let $\bk = (k_1,...,k_l)$, and let $\Delta_{\bk}$ be absolute value of the $\bk$-th Grassmann coordinate of $X$ and so of $\Lambda$ (i.e. $\Delta_{\bk}$ is absolute value of the $(N-l) \times (N-l)$ subdeterminant of $X$ obtained by removing the rows numbered $k_1,...,k_l$; this is an invariant of the lattice). Let $L_{k_1},...,L_{k_l}$ be distinct prime numbers so that $L_{k_i} \nmid \Delta_{\bk}$ for each $1 \leq i \leq l$. Define
\begin{equation}
\label{l1}
\Omega_{\bk} = \spn_{\zed} \{\bx_1,...,\bx_{N-1},L_{k_1}\be_{k_1},...,L_{k_l}\be_{k_l} \}.
\end{equation}
Then $\Lambda \subset \Omega_{\bk} \subseteq \zed^N$, $\Lambda \neq \Omega_{\bk}$, and $\Omega_{\bk}$ is a lattice of rank $N$. Notice that $\be_{k_1},...,\be_{k_l} \notin \Omega_{\bk}$.
\smallskip

Choose an integer basis $\ba_1,...,\ba_N$ for $\Omega_{\bk}$ so that the $N \times N$ basis matrix $A=(\ba_1 \hdots \ba_N)$ is upper triangular, and
\begin{equation}
\label{l2}
0 \leq a_{nj} < a_{nn}\ \ \forall\ 1 \leq n \leq N,\ 1 \leq j \leq N,\ j \neq n.
\end{equation}
Such a basis for $\Omega_{\bk}$ exists uniquely by Corollary 1 on p. 13 of \cite{cass:geom}. Notice that
$$\det(\Omega_{\bk}) = \Delta_{\bk} \prod_{i=1}^l L_{k_i} = |\det(A)| = \prod_{n=1}^N a_{nn}.$$
Fix an $s$, $1 \leq s \leq l$. Since $L_{k_s}$ is prime, $L_{k_s} | a_{nn}$ for some $1 \leq n \leq N$, and since $L_{k_s} \nmid \Delta_{\bk}$ this is the only $a_{nn}$ that $L_{k_s}$ divides.
\smallskip
Since $L_{k_s} \be_{k_s} \in \Omega_{\bk}$ and $A$ is upper triangular, there must exist integers $\alpha_{s1},...,\alpha_{sN}$ such that
$$L_{k_s} = \sum_{i=k_s}^N \alpha_{si} a_{ik_s},\ \ 0 = \sum_{i=j}^N \alpha_{si} a_{ij},\ \forall\ j\neq k_s.$$

{\it Case 1.} Suppose $k_s=N$. Then $L_N = \alpha_{sN} a_{NN}$, which implies that either $\alpha_{sN} = L_N,\ a_{NN} = 1$, or $\alpha_{sN} = 1,\ a_{NN} = L_N$. However, if $a_{NN}=1$, then by (\ref{l2}) $a_{Ni}=0$ for all $1 \leq i \leq N-1$, and so $\ba_N=\be_N \in \Omega_{\bk}$, which is a contradiction. Therefore $a_{NN}=L_N$.
\smallskip

{\it Case 2.} Suppose $k_s<N$. Then $\alpha_{sN} a_{NN}=0$, and so $\alpha_{sN}=0$. Then 
$$\alpha_{s(N-1)} a_{(N-1)(N-1)} + \alpha_{sN} a_{N(N-1)} = 0,$$
which means that $\alpha_{s(N-1)}=0$. Continuing in the same manner, we see that $\alpha_{si}=0$ for each $i>k_s$. Hence $L_{k_s} = \sum_{i=k_s}^N \alpha_{si} a_{ik_s} = \alpha_{sk} a_{k_sk_s}$. By the same argument as in case 1, this means that $a_{k_sk_s}=L_{k_s}$. 
\smallskip

Therefore we proved that $a_{k_sk_s}=L_{k_s}$ for all $1 \leq s \leq l$, and each $L_{k_s}$ does not divide any other $a_{nn}$, hence 
\begin{equation}
\label{l3}
\prod_{n=1,n \neq k_1,...,k_l}^N a_{nn} = \Delta_{\bk}.
\end{equation}
Applying Corollary \ref{2.4.1}, we see that for any $\bz \in \real^N$,
\begin{eqnarray}
\label{l4}
|\Lambda \cap (C_R^N+\bz)| & \leq & |\Omega_{\bk} \cap (C_R^N+\bz)| \nonumber \\
                           & \leq & \prod_{s=1}^l \left( \frac{2R}{L_{k_s}} + 1 \right) \prod_{n=1,n \neq k_1,...,k_l}^N \left( \frac{2R}{a_{nn}} + 1 \right).
\end{eqnarray}
Since our choice of $L_{k_1},...,L_{k_l}$ was arbitrary, we will now let $L_{k_s} \rightarrow \infty$ for all $1 \leq s \leq l$, and so
\begin{eqnarray}
\label{l5}
|\Lambda \cap (C_R^N+\bz)| & \leq & \prod_{n=1,n \neq k_1,...,k_l}^N \left( \frac{2R}{a_{nn}} + 1 \right) \prod_{s=1}^l \lim_{L_{k_s} \rightarrow \infty} \left( \frac{2R}{L_{k_s}} + 1 \right) \nonumber \\
                           & = & \prod_{n=1,n \neq k_1,...,k_l}^N \left( \frac{2R}{a_{nn}} + 1 \right).
\end{eqnarray}
The right hand side of (\ref{l5}) takes its maximum value when $a_{nn} = 1$ for $N-l-1$ distinct values of $n$. Therefore, applying (\ref{l3}) we obtain
\begin{equation}
\label{l6}
|\Lambda \cap (C_R^N+\bz)| \leq \left( \frac{2R}{\Delta_{\bk}} + 1 \right)(2R+1)^{N-l-1}.
\end{equation}
We can now specify how we select $\bk$. We want to do it so that the upper bound in (\ref{l6}) is minimized. For this, let $\bk$ be such that $\Delta_{\bk}$ is the maximal among all the Grassmann coordinates of $\Lambda$, and call this maximum value $\Delta$ (notice that if $\Delta_k \neq 0$, then the lattice $\spn_{\zed} \{\Lambda,\be_{k_1},...,\be_{k_l} \} \subseteq \zed^N$ has rank $N$). This completes the proof.  
\endproof

\noindent
The upper bounds in (\ref{cube:lattice}) and (\ref{notfull:rank}) are best possible as seen on the example of 
$$\Lambda = \spn_{\zed} \{ \be_1,...,\be_{N-l-1}, \Delta \be_{N-l} \},$$
where $0 \leq l \leq N-1$.
\smallskip

\begin{thm} \label{2.4.3} Suppose that $\Lambda \subseteq \zed^N$ is a lattice of rank $N-l$, where $1 \leq l \leq N-1$. Let $\Delta$ be the maximum of absolute values of Grassmann coordinates of $\Lambda$. Then for every $R$ that is a positive integer multiple of $(N-l) \Delta$, we have
\begin{equation}
\label{notfull:rank_lower}
\frac{(2R)^{N-l}}{(N-l)^{N-l} \Delta} \leq |\Lambda \cap C_R^N|.
\end{equation}
\end{thm}

\proof
Pick a basis $\bx_1,...,\bx_{N-l}$ for $\Lambda$ and write $X = (\bx_1\ ...\ \bx_{N-l})$ for the corresponding $N \times (N-l)$ basis matrix. Let $A$ be an $l \times N$ matrix such that $A\bx = \boldsymbol 0$ for each $\bx \in \Lambda$. Write $\Delta_{i_1...i_{N-l}}$ for the Grassmann coordinate of $X$, which is the determinant of the submatrix of $X$ whose rows are indexed by $i_1,...,i_{N-l} \in \{1,...,N\}$. Write $\delta^{j_1...j_l}$ for the Grassmann coordinate of $A$, which is the determinant of the submatrix of $A$ whose columns are indexed by $j_1,...,j_l \in \{1,...,N\}$. By the duality principle of Brill-Gordan \cite{gordan:1} (also see Theorem 1 on p. 294 of \cite{hodge:pedoe}) we have
\begin{equation}
\label{l7}
\Delta_{i_1...i_{N-l}} = (-1)^{i_1+...+i_{N-l}}\ \gamma\ \delta^{i_{N-l+1}...i_{N}},
\end{equation} 
for an appropriate $\gamma \in \real$, where $\{i_1,...,i_{N-l},i_{N-l+1},...,i_{N}\} = \{1,...,N\}$. We can assume without loss of generality that
\begin{equation}
\label{l8}
\Delta = |\Delta_{1...(N-l)}| = |\gamma| |\delta^{(N-l+1)...N}|.
\end{equation}
Let $X'$ be the $(N-l) \times (N-l)$ submatrix of $X$ whose rows are indexed by $1,...,(N-l)$, and let $\Lambda'$ be a lattice of full rank in $\real^{N-l}$ generated by the column vectors of $X'$. Then $\det(\Lambda') = \det(X') = \Delta$, and so, by the lower bound of Corollary \ref{2.4.1}
\begin{equation}
\label{l9}
|\Lambda' \cap C_R^{N-l}| \geq \frac{(2R)^{N-l}}{\Delta}.
\end{equation}
Suppose that $\bwy = (y_1,...,y_{N-l}) \in \Lambda' \cap C_R^{N-l}$, then $|\bwy| \leq R$, where $R \geq 1$. There exists $\bz = (z_1,...,z_l) \in \real^l$ such that $(\bwy, \bz) = (y_1,...,y_{N-l},z_1,...,z_l) \in \Lambda$. We want to establish an upper bound on $|\bz|$. By equation (4) on page 293 of \cite{hodge:pedoe}, every point $\bx \in \Lambda$ must satisfy the following system of linear equations:
\begin{equation}
\label{l10}
\sum_{j=1}^N \delta^{i_1...i_{l-1} j} x_j = 0,
\end{equation}
where $i_1,...,i_{l-1}$ assume all possible values; only $N-l$ of these equations are linearly independent. It is easy to see that the sum on the left side of each equation like (\ref{l10}) has only $N-l+1$ terms: there are only $N-l+1$ possibilities for $j$ since the $l-1$ values $i_1,...,i_{l-1}$ have been preassigned. For each $N-l+1 \leq i \leq N$ let $I_i = \{N-l+1,...,N\} \setminus \{i\}$, then the following $l$ equations form a subset of equations in (\ref{l10}):
\begin{equation}
\label{l11}
\sum_{j=1}^{N-l} \delta^{I_i j} x_j + \delta^{(N-l+1)...N} x_i = 0.
\end{equation}
Substitute the coordinates of the point $(\bwy,\bz)$ into (\ref{l11}), then we see that for each $1 \leq i \leq l$
\begin{equation}
\label{l12}
z_i = - \sum_{j=1}^{N-l} \frac{\delta^{I_{(N-l+i)} j}}{\delta^{(N-l+1)...N}} y_j,
\end{equation}
and so
\begin{equation}
\label{l13}
|z_i| \leq \sum_{j=1}^{N-l} \left| \frac{\delta^{I_{(N-l+i)} j}}{\delta^{(N-l+1)...N}} \right| |y_j| \leq (N-l) R,
\end{equation}
since each $\left| \frac{\delta^{I_{(N-l+i)} j}}{\delta^{(N-l+1)...N}} \right| \leq 1$ by construction, since $\delta^{(N-l+1)...N}$ is the biggest in absolute value among all the Grassmann coordinates of $A$, and $|\bwy| \leq R$. Therefore for each $\bwy \in \Lambda' \cap C_R^{N-l}$ there exists $\bz \in \real^l$ such that $(\bwy, \bz) \in \Lambda \cap C_{(N-l) R}^N$, hence
\begin{equation}
\label{l14}
|\Lambda \cap C_R^N| \geq |\Lambda' \cap C_{\frac{R}{N-l}}^{N-l}| \geq \frac{(2R)^{N-l}}{(N-l)^{N-l}\Delta}.
\end{equation}
This completes the proof.
\endproof


Notice that the entire argument in the proof of Theorem \ref{2.4.2}, except for (\ref{l9}), holds even if we drop the assumptions that $\Lambda \subseteq \zed^N$ and that $R$ is a positive integer multiple of $(N-l) \Delta$. In particular, this is true for (\ref{l14}).
\smallskip

{\it Remark.} The problem of counting lattice points in compact domains has been studied extensively. A long list of contributions on this subject can, for instance, be found in the supplement iv of Chapter 2 (pages 140 - 147) of \cite{lek}. The vast majority of such results is asymptotic in nature, which is not sufficient for our purposes. Explicit bounds on the number of lattice points in convex bodies can be found for instance in \cite{spain} and \cite{thunder}, however these bounds are, although quite general in the choice of the convex body, depend on parameters which are hard to compute. The advantage of bounds developed here in sections 3 and 4 is that they are reasonably sharp and easy to use in the particular case needed for our main result.
\bigskip

\section{Points outside of a collection of sublattices}

In this section we consider Problem \ref{problem2} and Problem \ref{problem3}. Throughout this section $N \geq 2$. Given a sublattice $\Omega$ of $\zed^N$ and a collection of proper sublattices $\Lambda_1,...,\Lambda_M$ of $\Omega$ we want to prove the existence of a non-zero integral point of small height in $\Omega \setminus \bigcup_{i=1}^M \Lambda_i$.
\smallskip 

We start with a discussion of Problem \ref{problem2}, namely the case when $\Omega=\zed^N$ and $\Lambda_1,...,\Lambda_M$ are sublattices of $\zed^N$ of rank $N-1$. There exist uniquely (up to $\pm$ sign) non-zero linear forms $L_1(\bX),...,L(\bX)$ in $N$ variables with integer relatively prime coeffcients such that
$$\Lambda_i = \{ \bwy \in \zed^N : L_i(\bwy) = 0 \}.$$
Then, given such a collection of linear forms, we want to prove the existence of an integer lattice point at which none of these linear forms would vanish.
\smallskip

A basic bound that depends only on the number of linear forms follows from results of section 2. Take $F(\bX) = L_1(\bX) \cdots L_M(\bX)$. By Lemma \ref{2.2.2}, there exists an integer lattice point $\bx$ such that $F(\bx) \neq 0$, and so $L_i(\bx) \neq 0$ for all $1 \leq i \leq N$, with
\begin{equation}
\label{basic1}
H(\bx) \leq \frac{M+1}{2}.
\end{equation}
Moreover, a bound of the form $H(\bx) \ll_N M^{(N-1)/N}$ follows from \cite{harcos}, however the constant there is not explicit.
\smallskip

We want to produce a result that depends on the actual linear forms, not just on their number. We will obtain such a result as a corollary of Theorem \ref{classical:subspace}, which presents a solution to the more general Problem \ref{problem3}. Namely, we prove the existence of a point of small height in a sublattice of the integer lattice outside of a collection of its proper sublattices. We will now relate this problem to the lattice point counting problem of sections~3 and~4 to prove our main theorem.
\bigskip

{\it Proof of Theorem \ref{classical:subspace}.}
For each $1 \leq i \leq M$ let $l_i$ be the rank of lattice $\Lambda_i$, then $0 < l_i < J$. Define a counting function of a variable $R$, which is a positive integer multiple of $J H(\Omega)$, by
\begin{equation}
\label{c4}
f_{\Omega}(R) = |\Omega \cap C_R^N| - \left| \bigcup_{i=1}^M \Lambda_i \cap C_R^N \right|,
\end{equation}
then
\begin{eqnarray}
\label{c5}
f_{\Omega}(R) & \geq & |\Omega \cap C_R^N| - \sum_{i=1}^M |\Lambda_i \cap C_R^N| \nonumber \\
       & \geq & \frac{(2R)^J}{J^J H(\Omega)} - \sum_{i=1}^M \left( \frac{2R}{H(\Lambda_i)}+1 \right)(2R+1)^{l_i-1},
\end{eqnarray}
where the last inequality follows by Corollary \ref{2.4.1}, Theorem \ref{2.4.2}, and Theorem \ref{2.4.3}. Notice that $f_{\Omega}(R)>0$ if and only if there exists a point $\bx \in \Omega \setminus \bigcup_{i=1}^M \Lambda_i$ with $H(\bx) \leq R$. Hence ideally we want to find $R$, the smallest positive integer multiple of $J H(\Omega)$, so that $f_{\Omega}(R)>0$. Recall that $l_i \leq J-1$ for each $i$. Then using (\ref{c5}), we have
$$f_{\Omega}(R) \geq R^{J-2} \left\{ \left(\frac{2^J}{J^J H(\Omega)}\right) R^2 - 3^{J-1} \left(\sum_{i=1}^M \frac{1}{H(\Lambda_i)}\right) R - 3^{J-2} M \right\}.$$
Therefore we want to solve for $R$ a quadratic inequality
\begin{equation}
\label{c6}
\left(\frac{2^J}{J^J H(\Omega)}\right) R^2 - 3^{J-1} \left(\sum_{i=1}^M \frac{1}{H(\Lambda_i)}\right) R - 3^{J-2} M > 0.
\end{equation}
It is not difficult to deduce from (\ref{c6}) that we can take $R$ to be a positive integer multiple of $J H(\Omega)$ such that
\begin{equation}
\label{c7}
R \geq \left( \frac{3}{2} \right)^{J-1} J^J \left\{ \sum_{i=1}^M \frac{1}{H(\Lambda_i)} + \sqrt{M} \right\} H(\Omega).
\end{equation}
This completes the proof.
\boxed{ }
\bigskip

The dependence on $H(\Omega)$, $H(\Lambda_i)$ for each $i$, and on $M$ in the upper bound of Theorem \ref{classical:subspace} seems to be best possible. Let $M=1$, and take $\Lambda_1$ to be a sublattice of $\Omega$ of rank $J-1$ generated by the vectors corresponding to the first $J-1$ successive minima of $\Omega$ with respect to a unit cube. Then the smallest vector in $\Omega \setminus \Lambda_1$ will be the one corresponding to the $J$-th successive minimum, and its height can be approximated by $H(\Omega)$. The dependence on $H(\Lambda_i)$ is sharp because it comes from the upper bounds on the number of lattice points in a cube (see (\ref{cube:lattice}) and (\ref{notfull:rank})), which are best possible. Finally, the proof of Theorem \ref{classical:subspace} (specifically see (\ref{c6})) suggests that dependence on $M$ is essentially sharp.
\bigskip

 


Also notice that since Problem \ref{problem2} is a special case of Problem \ref{problem3}, a solution to it follows from Theorem \ref{classical:subspace} when $J=N$, i.e. $\Omega = \zed^N$, and $l_i = N-1$ for all $1 \leq i \leq M$. In this case we obtain a point $\bx \in \zed^N$ such that $\bx \notin \bigcup_{i=1}^M \Lambda_i$ and
$$H(\bx) \leq \left( \frac{3}{2} \right)^{N-1} N^N \left\{ \sum_{i=1}^M \frac{1}{H(\Lambda_i)} + \sqrt{M} \right\} + N.$$
Hence Theorem \ref{classical:subspace} can be thought of as a combination of Siegel's Lemma and its inverse problem of finding points of small height outside of a collection of sublattices. 
\smallskip

\noindent
In fact, a slightly better bound for the point in question in Problem \ref{problem2} can be obtained from the proof of Theorem \ref{classical:subspace}. Notice that if $\Omega=\zed^N$ then in the lower bound of (\ref{c5}) the quantity $\frac{(2R)^J}{J^J H(\Omega)} = \left( \frac{2R}{N} \right)^N$ can be replaced with the slightly larger $(2R+1)^N$, since $|\zed^N \cap C_R^N| = (2R+1)^N$ for positive integer values of $R$. This leads to the following slightly sharper bound. 

\begin{cor} \label{classical:main} Let $L_1(\bX),...,L_M(\bX) \in \zed[X_1,...,X_N]$ be non-zero linear forms with relatively prime coordinates. Then there exists $\bx \in \zed^N$ such that $L_i(\bx) \neq 0$ for every $i=1,...,M$ and
\begin{equation}
\label{classical:bound}
H(\bx) \leq \sum_{i=1}^M \frac{1}{H(L_i)} + \sqrt{M}.
\end{equation}
\end{cor}

Corollary \ref{classical:main} produces a better bound than (\ref{basic1}) for linear forms with suffciently large heights. Also, suppose that $N$ is fixed and $M$ grows. Then our collection must contain linear forms with relatively large heights, since there are only finitely many vectors of height $\leq C$ in $\zed^N$ for each $C$. This is definitely the more interesting situation, since if $M<N$ or if the two are comparable, there must exist integer lattice points of height $\ll N$ at which the linear forms do not vanish.
\bigskip

Another immediate corollary of Theorem \ref{classical:subspace} in case when $M=1$ is a sharper bound for a special case of Faltings' version of Siegel's Lemma (see \cite{faltings}, \cite{faltings:siegel}, and \cite{faltings:siegel_1}). 

\begin{cor} \label{faltings:lemma} Let $V$ and $W$ be real vector spaces of respective dimensions $d_1$ and $d_2$. Let $\Omega_1 = V \cap \zed^{d_1}$ and $\Omega_2 = W \cap \zed^{d_2}$. Let $\rho: V \longrightarrow W$ be a linear map such that $\rho(\Omega_1) \subseteq \Omega_2$. Let $U = \ker(\rho)$, and let $\Omega = U \cap \Omega_1$. Let $J$ be the rank of $\Omega$. Then for any subspace $U_0$ of $U$ which does not contain $\Omega$ there exists a point $\bx \in \Omega \setminus U_0$ such that
\begin{equation}
\label{faltings:bound_my}
H(\bx) \leq \left\{ 2 J^J \left( \frac{3}{2} \right)^{J-1} + J \right\} H(\Omega).
\end{equation}
\end{cor}

\noindent
Faltings' lemma is more general: it works with $\Omega_1$, $\Omega_2$, and $\Omega$ being any lattices in $V$, $W$, and $U$ respectively, as well as any choice of norms on $V$ and $W$ (our height is the sup-norm). However, the upper bound on $H(\bx)$ which follows from Faltings' lemma is 
\begin{equation}
\label{faltings:bound}
H(\bx) \leq d_1^{\frac{1}{d - d_0}} H(\Omega)^{\frac{3d_1}{d - d_0}},
\end{equation}
where $d=\dim_{\real}(U)$ and $d_0=\dim_{\real}(U_0)$, so that $d_1 > d > d_0$ and $d \geq J$. Faltings' method of proof is different from ours: it relies on Minkowski's theorem about successive minima.
\bigskip

\section{Tarski plank problem}

We now consider a simple application of our results to a certain analogue of the discrete version of the Tarski plank problem. First we provide some background. By a {\it plank} of width $h$ in $\real^N$ we mean a strip of space of width $h$ between two parallel $(N-1)$-dimensional hyperplanes. Let $C$ be a convex body of minimal width $w$ in $\real^N$. If $C$ is covered by $p$ planks of widths $h_1,...,h_p$ respectively, is it true that $h_1+...+h_p \geq w$? This question was originally asked by Tarski in \cite{tarski}. It was answered affirmitively by Bang in \cite{bang}. One can also ask for the minimal number of planks with prescribed widths that would cover $C$. The discrete version of this problem (see for instance \cite{corzatt}) asks for the minimal number of $(N-1)$-dimensional hyperplanes that would cover a convex set of lattice points in $\real^N$. We ask a somewhat different, but analogous question. Consider the set of all integer lattice points in $\real^N$ that are contained in the closed cube $C_R^N$, where $R$ is a positive integer as above. This set has cardinality $(2R+1)^N$. What is the minimal number of sublattices of $\zed^N$ of rank $N-1$ required to cover this set? Let $M$ be this number. Then the inequality
\begin{equation}
\label{tarski1}
M \geq 2R-1
\end{equation}
follows immediately from (\ref{basic1}). For such a sublattice $\Lambda$ let $L(\bX) \in \zed[X_1,...,X_N]$ be the linear form with relatively prime coefficients such that $\Lambda=\{\bwy \in \zed^N : L(\bwy)=0\}$. An analogue of width of a plank in this case would be the quantity $H(L)^{-1} = H(\Lambda)^{-1}$, and the sidelength of the cube $C_R^N$ which is equal to $2R$ is an analogue of the width of a convex body. Then we can state the following result, which is an immediate corollary of Corollary \ref{classical:main}.

\begin{cor} \label{tarski:cor} Let $\Lambda_1,...,\Lambda_M$ be sublattices of $\zed^N$ of rank $N-1$ each, such that $C_R^N \cap \zed^N \subset \bigcup_{i=1}^M \Lambda_i$. Then
\begin{equation}
\label{tarski:bound}
\sum_{i=1}^M \frac{1}{H(\Lambda_i)} \geq R-\sqrt{M}.
\end{equation}
\end{cor}
\smallskip

A similar problem is treated in \cite{harcos}. Let $C$ be a compact convex body, which is symmetric with respect to the origin in $\real^N$. Suppose that $C$ can be inscribed into a cube $C_R^N$ as above. How many $(N-1)$-dimensional subspaces of $\real^N$ does it take to cover $C \cap \zed^N$, the set of integer lattice points contained in $C$? Call this number $M$. Theorem 2 of \cite{harcos} provides an upper bound for $M$ in terms of the successive minima of $C$ with respect to $\zed^N$, which implies that $M$ is of the order of magnitude $O(R^{N/(N-1)})$, which is better than (\ref{tarski1}). However, the actual constants in the inequalities of \cite{harcos} are not effectively computable, since they rely on successive minima. More precisely, Theorem 2 of \cite{harcos} states that
\begin{equation}
\label{harcos:bound}
M \leq c 2^N N^2 \log N \min_{0 < m < N} (\lambda_m \cdots \lambda_N)^{-\frac{1}{N-m}},
\end{equation}
where $0 < \lambda_1 \leq ... \leq \lambda_N \leq 1$ are the successive minima of $C$ with respect to $\zed^N$ (the case $\lambda_N > 1$ is trivial: $M=1$), and $c$ is an absolute constant. 
\bigskip



{\bf Aknowledgements.} I would like to express my deep gratitude to Professor Jeffrey D. Vaaler for his valuable advice and numerous helpful conversations on the subject of this paper. In particular, the argument of section 3 was suggested by him. I would also like to thank Dr. Iskander Aliev and Professor Preda Mihailescu for their very helpful remarks.

\nocite{*}
\bibliographystyle{plain}  
\bibliography{fukshansky}        

\begin{thebibliography}{10}

\bibitem{bang}
T.~Bang.
\newblock A solution of the ``{P}lank {P}roblem''.
\newblock {\em Amer. Math. Soc. Proceedings}, 2:990--993, 1951.

\bibitem{harcos}
I.~B\'ar\'any, G.~Harcos, J.~Pach, and G.~Tardos.
\newblock Covering lattice points by subspaces.
\newblock {\em Period. Math. Hungar.}, 43:93--103, 2001.

\bibitem{vaaler:siegel}
E.~Bombieri and J.~D. Vaaler.
\newblock On {S}iegel's lemma.
\newblock {\em Invent. Math.}, 73(1):11--32, 1983.

\bibitem{local:riemann}
D.~Bump, K.~K. Choi, P.~Kurlberg, and J.~D. Vaaler.
\newblock A local {R}iemann hypothesis. {I}.
\newblock {\em Math. Z.}, 233(1):1--19, 2000.

\bibitem{cass:geom}
J.~W.~S. Cassels.
\newblock {\em An Introduction to the Geometry of Numbers}.
\newblock Springer-Verlag, 1959.

\bibitem{corzatt}
C.~Corzatt.
\newblock Covering convex sets of lattice points with straight lines.
\newblock {\em Proceedings of the Sundance conference on combinatorics and
  related topics. Congr. Numer.}, 150:129--135, 1985.

\bibitem{faltings:siegel_1}
B.~Edixhoven.
\newblock Arithmetic part of {F}altings's proof.
\newblock {\em Diophantine approximation and abelian varieties (Soesterberg,
  1992)}, Lecture Notes in Math.(1566):97--110, 1993.

\bibitem{faltings}
G.~Faltings.
\newblock Diophantine approximation on abelian varieties.
\newblock {\em Ann. of Math.}, 133(2):549--576, 1991.

\bibitem{me:number}
L.~Fukshansky.
\newblock Siegel's lemma with additional conditions.
\newblock {\em submitted for publication}.

\bibitem{me:smallzeros}
L.~Fukshansky.
\newblock Small zeros of quadratic forms with linear conditions.
\newblock {\em to appear in J. Number Theory}.

\bibitem{me:diss}
L.~Fukshansky.
\newblock {\em Algebraic points of small height with additional arithmetic
  conditions}.
\newblock PhD thesis, {U}niversity of {T}exas at {A}ustin, 2004.

\bibitem{gordan:1}
P.~Gordan.
\newblock Uber den grossten gemeinsamen factor.
\newblock {\em Math. Ann.}, 7:443--448, 1873.

\bibitem{lek}
P.~M. Gruber and C.~G. Lekkerkerker.
\newblock {\em Geometry of Numbers}.
\newblock North-Holland Publishing Co., 1987.

\bibitem{hodge:pedoe}
W.~V.~D. Hodge and D.~Pedoe.
\newblock {\em Methods of Algebraic Geometry, Volume 1}.
\newblock Cambridge Univ. Press, 1947.

\bibitem{faltings:siegel}
R.~J. Kooman.
\newblock Faltings's version of {S}iegel's lemma.
\newblock {\em Diophantine approximation and abelian varieties (Soesterberg,
  1992)}, Lecture Notes in Math.(1566):93--96, 1993.

\bibitem{siegel}
C.~L. Siegel.
\newblock Uber einige {A}nwendungen diophantischer {A}pproximationen.
\newblock {\em Abh. der Preuss. Akad. der Wissenschaften Phys.-math Kl.}, Nr.
  1:209--266, 1929.

\bibitem{spain}
P.~G. Spain.
\newblock Lipschitz: a new version of an old principle.
\newblock {\em Bull. London Math. Soc.}, 27:565--566, 1995.

\bibitem{tarski}
A.~Tarski.
\newblock Further remarks about the degree of equivalence of polygons (in
  {P}olish).
\newblock {\em Odbitka Z. Parametru.}, 2:310--314, 1932.

\bibitem{thue}
A.~Thue.
\newblock Uber {A}nnaherungswerte algebraischer {Z}ahlen.
\newblock {\em J. Reine Angew. Math.}, 135:284--305, 1909.

\bibitem{thunder}
J.~L. Thunder.
\newblock The number of solutions of bounded height to a system of linear
  equations.
\newblock {\em J. Number Theory}, 43:228--250, 1993.

\end{thebibliography}

\end{document}